\documentstyle[11pt,twoside]{article}
\newfont{\BBB}{msbm10 at 18pt}
\newfont{\BBb}{msbm10 at 14pt}
\newfont{\Bbb}{msbm10 at 11pt}
\newfont{\Bbs}{msbm10 at 8pt}

\newtheorem{The}{Theorem}[section]
\newtheorem{Def}[The]{Definition}
\newtheorem{Pro}[The]{Proposition}

\newtheorem{Cor}[The]{Corrolary}

\begin{document}
\pagestyle{myheadings}
\markboth{M. PAPATRIANTAFILLOU}{HERMITIAN CONNECTIONS ON
            {\Bbb A}-BUNDLES}

\title{HERMITIAN STRUCTURES AND COMPATIBLE CONNECTIONS ON
                        {\BBB A}-BUNDLES$^*$}
\author{M. Papatriantafillou}
\date{}
\maketitle
\thispagestyle{empty}

\noindent
{\bf AMS subject classification:} 53C05, 58B20

\noindent
{\bf Key words:} {\Bbb A}-bundles, hermitian structures, connections

\medskip
\noindent
{\bf *} This paper is in final form and no version of it will be
           submitted for publication elsewhere

\def\A{\mbox{\Bbb A}}
\def\a8{{\Bbb A}^{\infty}}
\def\c{\mbox{\Bbb C}}
\def\C{{\cal C}^{\infty}}
\def\G{\Gamma(l)}
\def\e{\varepsilon}
\def\R{\mbox{\Bbb R}}
\def\N{\mbox{\Bbb N}}
\def\n{{\cal N}}
\def\P{{\cal P}(\A)}
\def\ra{\rightarrow}
\def\mt{\mapsto}
\def\eeq{\end{equation}}
\def\beq{\begin{equation}}

\section{Introduction}

Arbitrary locally convex spaces have a very poor geometric structure
and this fact is reflected to the manifolds and vector bundles
modelled on them (cf., for instance, \cite{Michor,Papaghiuc}). In
particular, they do not have inner products, thus a  vector bundle
modelled on a locally convex space is not endowed with a Riemannian
structure.

But some locally convex spaces, arising in pure mathematics and in
theoretical physics, have the additional algebraic structure of a
(projective finitely generated) module over a topological algebra
$\A$ and a number of questions have been answered by the extension of
the usual ring $\R$ or $\c$ of coefficients to the aforementioned
algebra (see \cite{KAPL} in operator theory, \cite{SEL} in
theoretical physics, \cite{MISC-SOL} in differential topology).
In the case that $\A$ is a $^*$-algebra, the modules are provided
with $\A$-valued inner products and norms, defining their topology
\cite{MP;Buda}, and they behave like finite dimensional vector
spaces, although, in general, they lack both bases and metric
topologies.

On the other hand, manifolds and vector bundles modelled on such
modules are found in various areas (see, for example, \cite{MISC-SOL}
in differential topology, \cite{PRAST} in partial differential
equations, \cite{KOB} in global analysis, \cite{SHUR} in the theory of
jets), $\A$ usually being ${\cal C}(X)$ or $\C(X)$. For brevity, we
call them $\A$-manifolds and $\A$-bundles, respectively.

Continuous $\A$-bundles have been extensively studied (see
\cite{MAL1,MAL2,MAL3,MAL4,MP;SS,MP;Red}), while some differential
aspects have appeared in \cite{MP;Finsler,MP;JMAA,MP;Buda}, among
which the existence of $\A$-valued Finsler structures. Our aim here is
to investigate the conditions under which an $\A$-bundle is provided
with generalized (: $\A$-valued) hermitian structures and compatible
connections, in the general case when $\A$ is a {\em commutative
locally m-convex} $^*$-{\em algebra with unit}. In this investigation
two obstacles appear: first, $\A$-manifolds never admit partitions of
unity in the classical sense, and, secondly, the existence of a
hermitian structure is not equivalent to the reduction of the
structural group of the bundle to a special subgroup. However, we
prove that an $\A$-hermitian structure exists if the fibre type of
the $\A$-bundle has an $\A$-hermitian inner product and the base
space admits just one $\A$-valued partition of unity (Theorem 4.5),
or, if the structural group of the bundle reduces to the
``$\A$-hermitian product preserving'' automorphisms (Theorem 4.7).
Next, we endow an $\A$-bundle with a connection, assuming the
existence of one $\A$-partition of unity (Theorem 5.3), and we prove
that this connection and the $\A$-hermitian structure of Theorem 4.7
are compatible (Theorem 5.4).

\section{Preliminaries}

We recall that a complex algebra $\A$ is a $^*$-{\em algebra},
if it is endowed with a map $^* : \A \to \A$, so that (i) $(za+b)^* =
\bar{z}a^* + b^*$, (ii) $(ab)^* = b^*a^*$ and (iii) $(a^*)^* = a$, for
every $a,b \in \A$, $z \in \c$. A $^*$-algebra $\A$ is called a
{\em locally m-convex} (abr. {\em lmc}) $^*$-{\em algebra}, if it is
topologized by a family of seminorms that satisfy (i) $p(xy)
\leq p(x)p(y)$ and (ii) $p(x^*) = p(x)$, for every $x,y \in \A$ (for
details, see \cite{MAL;TA}). 

Throughout the paper, $\A$ denotes a {\em commutative lmc\/})
$^*$-{\em algebra with unit\/}.

Let $\P$ be the category of projective finitely generated
$\A$-modules. By definition, for any $M \in \P$, there exist
$M_1 \in \P$ and $m \in \N$, so that $M \oplus M_1 \cong \A^m$. Let
$\tau_M$ denote the {\em canonical topology} of $M$, i.e., the
relative topology induced on $M$ by the product topology of $\A^m$.
Then: (i) $\tau_M$ is independent of either $M_1$ or $m$; (ii)
$(M,\tau_M)$ is a topological $\A$-module (namely, the $\A$-module
operations are jointly continuous); (iii) $\tau_M$ makes every
$\A$-linear map $f : M \ra N$ continuous, for any topological
$\A$-module $N$ (for details, see \cite{MP;PMH}).

In the sequel, every $M \in \P$ is a topological $\A$-module
provided with the canonical topology. For any $M, \, N \in \P$ and
$x \in M$, we denote by $0_M$ the zero element of $M$, by $\n(x)$
the set of open neighbourhoods of $x$ and by $L_{\mbox{\Bbs A}}(M,N)$
(resp. $S_{\mbox{\Bbs A}}(M,N)$) the set of $\A$-linear (resp.
skew-linear) maps $f : M \ra N$; we recall that a map
$f : M \ra N$ is called {\em skew-linear}, if $f$ is additive
and $f(ax) = a^*f(x)$, for every $x \in M$ and $a \in \A$.

\medskip
In $\P$ we consider the following differentiation method: Let
$M, \, N \in \P$, $x \in M$, $U \in {\cal N}(x)$ and $f : U \ra N$.
We say that $f$ is $\A$-{\em differentiable at} $x$, if there exist
$Lf(x) \in L_{\mbox{\Bbs A}}(M,N)$ and $Sf(x) \in
S_{\mbox{\Bbs A}}(M,N)$, such that the map
\[
\phi(h) := f(x+h) - f(x) - Lf(x)(h) - Sf(x)(h)
\]
satisfies the following condition:
\[
\forall \, V \in \n(0_N) \; \exists \, U \in \n(0_M) \, : \,
\forall \, B \in \n(0_{\mbox{\Bbs A}}) \; \exists \,
A \in \n(0_{\mbox{\Bbs A}}) :
\]
\[
\phi(aU+a^*U) \subseteq aBV+a^*BV, \quad \forall a \in \A.
\]
We call $Df(x) := Lf(x)+Sf(x)$ {\em the differential of} $f$ at $x$.
If $Sf(x) = 0$ (resp. $Lf(x) = 0$), $f$ is called
$\A$-{\em holomorphic} (resp. $\A$-{\em antiholomorphic}) {\em at}
$x$.

Let $f$ be $\A$-differentiable at every $x \in U$. Since
$L_{\mbox{\Bbs A}}(M,N) \oplus S_{\mbox{\Bbs A}}(M,N) \in \P$,
$\A$-differentiation may apply to
\[
Df : U \ra L_{\mbox{\Bbs A}}(M,N) \oplus S_{\mbox{\Bbs A}}(M,N),
\]
inducing {\em the second differential} $D^2f = D(Df)$ {\em of} $f$,
and, successively, {\em the} $n${\em -th differential} $D^nf$
{\em of} $f$, for any $n \in \N$. We will say that $f$ is an
$\A^{\infty}$-{\em differentiable} (resp.
$\A^{\infty}$-{\em holomorphic}) {\em map on} $U$, if $D^kf$ exists
(resp. $D^kf$ exists and $S^kf = 0$), for every $k \in \N$ (for
details, we refer the reader to \cite{MP;JMAA}).

\section{{\BBb A}-manifolds and {\BBb A}-bundles}

Let $X$ be a Hausdorff topological manifold modelled on $M \in \P$. We
say that $X$ is an $\A$-{\em manifold} (resp. $\A_h$-{\em manifold}),
if its transition functions are $\A^{\infty}$-differentiable (resp.
$\a8$-holomorphic). If $X, \, Y$ are $\A$-manifolds (resp.
$\A_h$-manifolds), we say that $f : X \ra Y$ is an $\A$-{\em map}
(resp. $\A_h$-{\em map}), if its local representatives are
$\A^{\infty}$-differentiable (resp. $\a8$-holomorphic). The category
of $\A$-manifolds (resp. $\A_h$-manifolds) and $\A$-maps (resp.
$\A_h$-maps) will be denoted by $Man(\A)$ (resp. $Man_h(\A))$. 

Let $X \in Man_h(\A)$ modelled on $M$. We obtain tangent spaces, by
considering classes of equivalent ``curves'' in the following way:
an $\A$-{\em curve on} $X$ is an $\A$-map $\alpha : A \ra X$, with
$A \in \n(0_{\mbox{\Bbs A}})$. The $\A$-curves $\alpha, \, \beta$
{\em are tangent at} $x \in X$, if $\alpha(0) = \beta(0) = x$ and
there exists a chart $(U,\phi)$ at $x$ with
$D(\phi \circ \alpha)(0) = D(\phi \circ \beta)(0)$. We denote by
$[(\alpha,x)]$ the induced equivalence class of $\alpha$ and by
$T(X,x)$ the set of such quivalence classes. If $M_*$ denotes the
abelian group $(M,+)$ provided with the scalar multiplication
\[
\A  \times  M_* \ra M_* : (a,x) \ra a^*x
\]
and $M \oplus M_1 = \A^m$, then
\[
M_* \oplus (M_1)_* = (M \oplus M_1)_* = (\A^m)_* = \A^m,
\]
within $\A$-module isomorphisms, that is, $M_* \in \P$. Let $x \in X$
and $(U,\phi)$ a chart at $x$. The map
\begin{equation}
\bar{\phi} : T(X,x) \ra M  \times  M_* : [(\alpha,x)] \ra
(L(\phi \circ \alpha)(0),S(\phi \circ \alpha)(0))
\end{equation}
is a bijection establishing an $\A$-module structure on $T(X,x)$.
We call $T(X,x)$ {\em the tangent space of} $X$ {\em at} $x$. The
{\em tangent bundle} $T(X)$ {\em of} $X$, i.e., the discrete union of
all tangent spaces is an $\A$-manifold.

We note here, that if $\A = \c$, the tangent bundle introduced above
coincides with the complexified tangent bundle of complex manifolds.

If  $f : X \ra Y$  is an $\A$-map, the {\em differential of} $f$
\[
Tf : T(X) \ra T(Y) : [(\alpha,x)] \ra [(f \circ \alpha,f(x))]
\]
is an $\A$-map  and, for any $x \in X$, the restriction
\[
T_xf : T(X,x) \ra T(Y,f(x)) : [(\alpha,x)] \ra [(f \circ \alpha,f(x))]
\]
is an $\A$-linear map.

\medskip
Let now $X \in Man_h(\A)$,  $E \in Man(\A)$, $\pi : E \ra  X$ an
$\A$-map  and $M \in \P$. We say that the triplet $\ell = (E,\pi,X)$
is an $\A$-{\em differentiable} $\A$-{\em bundle over $X$ of fibre
type} $M$, or, simply, an $\A$-{\em bundle}, if the following
conditions hold:

i) $E_x := \pi^{-1}(x) \in \P$, for every $x \in X$.

ii) There exists a {\em trivializing covering}
$\{ (U_i,\tau_i) \}_{i \in I}$, where $\{U_i\}_{i \in I}$  is an open
covering of $X$ and every $\tau_i : \pi^{-1}(U_i) \ra U_i  \times  M$
is an isomorphism in $Man(\A)$, such that $pr_1 \circ \tau_i = \pi$
and, for every $x \in U_i$, the restriction
\[
\tau_{ix} := \tau_i|_{E_x} : E_x \ra \{ x \}  \times  M
\]
is an $\A$-module isomorphism.

\medskip
One would note here that in the Banach context, for the definition of
vector bundles, one more condition is required, namely (VB 3) of
\cite{LANG}. However, in our framework, the properties of the
canonical topology and the $\a8$-differentiation imply this
condition, making $\A$-bundles look like bundles of finite rank (cf.
the analogous results for continuous $R$-bundles, where $R$ is a
topological ring \cite{MP;SS} and for differentiable $\A$-bundles,
where $\A$ is a commutative unital lmc algebra over $\R$
\cite{MP;Buda,MP-EV}).

\section{{\BBb A}-hermitian structures}

The involution of the algebra $\A$ endows the objects of $\P$ with a
structure generalizing hermitian inner products on  complex vector
spaces. In this section we investigate the conditions under which
these generalized inner products provide a hermitian structure on an
$\A$-bundle.

\medskip
We recall that a {\em positive definite} $\A$-{\em hermitian inner
product} on an $\A$-module $M$ is a map
$\alpha : M  \times  M \ra \A$, satisfying the following conditions:

(i) $\alpha$ is $\A$-linear with respect to the first variable. 

(ii) $\alpha(y,x) = (\alpha(x,y))^*$, for any $x,y \in M$.

(iii) For every $x \in M$, $\alpha(x,x)$ is positive in $\A$, that is,
\[
sp_{\mbox{\Bbs A}}(\alpha(x,x)) :=
\{ \lambda \in {\mbox{\Bbb C}} : \lambda \cdot 1_{\mbox{\Bbs A}} -
\alpha(x,x) \mbox{ not invertible } \} \subseteq [0,+\infty).
\]

(iv) The mapping $M \ra L_{\mbox{\Bbs A}}(M,\A)_* : x \ra \alpha_x$,
where $\alpha_x(y) := \alpha(y,x)$, for every $y \in M$, is an
isomorphism of $\A$-modules.

For brevity, the pair $(M,\alpha)$ is called a {\em hermitian form}.

\medskip
For every hermitian form $(M,\alpha)$, with $M \in \P$, $\alpha$ is an 
$\A$-map and 
\[
L\alpha(x,y)(h,k) = \alpha(h,y), \qquad
S\alpha(x,y)(h,k) = \alpha(x,k),
\]
for every $(x,y), (h,k) \in M  \times  M$.

\medskip
Let us recall that a {\em lmc C}$^*$-{\em algebra} is a lmc
$^*$-algebra, whose seminorms satisfy the relation $p(x^*x) =
(p(x))^2$, for every $x \in \A$. Regarding the existence of
hermitian forms, we have

\begin{The}  \label{exist}
{\em \cite{MAL2}} Let $\A$  be a complete lmc $C^*$-algebra with unit
and $M \in \P$. Then $M$  admits a positive definite $\A$-hermitian
inner product $\alpha$,  which is unique up to an isomorphism, that
is, if $(M,\beta)$ is a hermitian form, then there exists an
$\A$-automorphism $f$ of $M$, such that $\beta \circ (f \times f) =
\alpha$. \hfill $\Box$
\end{The}

\begin{Cor}
Let $\A$  be a complete lmc $C^*$-algebra with unit and $M$ a free
finitely generated $\A$-module. If $(M,\beta)$ is a hermitian form, 
Then $M$ has an orthonormal basis with respect to $\beta$.
\end{Cor}

\noindent {\em Proof.} 
By definition, $M$ coincides with $\A^m$, for some $m \in \N$. The
canonical basis $\{e_i\}_{i=1,\dots,m}$ of $\A^M$ is orthonormal with
respect to the posotive definite $\A$-hermitian inner product
\[
\alpha : \A^m \times \A^m \to \A : ((x_i),(y_i)) \mapsto
\sum_i x_i(y_i)^*.
\]
If $f$ is the $\A$-automorphism of $\A^m$ with $\beta \circ
(f \times f) = \alpha$, then $\{f(e_i)\}_{i=1,\dots,m}$ is the
required basis of $\A^m$. \hfill $\Box$

\begin{Def}
{\em An} $\A$-hermitian structure on the $\A$-bundle $\ell =
(E,\pi,X)$ {\em is an} $\A$-{\em map}
$g : E \oplus E \ra \A${\em , such that,}
$(E_x,g_x := g|_{E_x \oplus E_x})$ {\em is a hermitian form,
for every} $x \in X$.
\end{Def}

In the subsequent theorems we give sufficient conditions for the 
existence of $\A$-hermitian structures. But, we need first the 
following

\begin{Def}
{\em Let} $(X,{\cal A}) \in Man(\A)$. {\em An} $\A$-partition of unity
on $X$ {\em is a family} $\{ (U_i,\psi_i) \}_{i \in I}$, {\em where}
$\{ U_i \}_{i \in I}$ {\em is a locally finite open covering of}
$X$ {\em and} $\{ \psi_i \}_{i \in I}$ {\em is a family of}
$\A$-{\em maps} $\psi_i : X \ra \A$, {\em such that

(i)} $supp(\psi_i) \subseteq U_i$, {\em for any} $i \in I$.

{\em (ii) For any} $x \in X$ {\em and any} $i \in I$, $\psi_i(x)$ {\em
is positive in} $\A$.

{\em (iii)} $\sum_i \psi_i(x) = 1$, {\em for any} $x \in X$.
\end{Def}

In the ordinary (finite dimensional or Banach) case, one assumes that
every open covering of $X$ admits a locally finite refinement with
a subordinate partition of unity. However, such an assumption is too 
strong. We only need that the bundle has (at least) one locally finite 
trivializing covering $\{ (U_i,\tau_i) \}_{i \in I}$, so that
$\{ U_i \}_{i \in I}$ has a subordinate $\A$-partition of unity. This
last condition is proved to hold for a class of $\A$-manifolds, if
$\A$ is the algebra ${\cal C}(X)$ of continuous complex valued
functions on a Hausdorff completely regular topological space $X$,
or its subalgebra $\C(X)$ of smooth functions, in the case that $X$
is a compact smooth manifold (see \cite{MP;Bumpkob}, \cite{MP;Bump}).

Following the classical arguments, one has

\begin{The}
Let $\ell$ be an $\A$-bundle of fibre type $M \in \P$. If
$(M,\alpha)$ is an $\A$-hermitian form and $\ell$ has a localy finite
trivializing covering with a subordinate $\A$-partition of unity, then
$\ell$ is provided with an $\A$-hermitian structure. \hfill $\Box$
\end{The}

In the case of a unital commutative complete lmc C$^*$-algebra, for
every $M \in \P$ there is a hermitian form $(M,\alpha)$ (Theorem
\ref{exist}). Besides, if the base space $X$ has a locally finite
atlas consisting of charts whose image is a sphere with respect to
the $\A$-valued norm defined by $\alpha$, then $X$ has an
$\A$-partition of unity subordinate to this atlas (see \cite{MP;Bump}
in conjunction with \cite{MP;Bumpkob}). As a result, we obtain

\begin{The}
Let $\A$ be a unital commutative complete lmc C$^*$-algebra. If
the $\A$-bundle $\ell = (E,\pi,X)$ has a localy finite trivializing
covering, so that the respective charts of $X$ are sent to spheres,
then $\ell$ is provided with an $\A$-hermitian structure.
\hfill $\Box$
\end{The}

If $(M,\alpha)$ is a hermitian form, we denote by $GL(M,\alpha)$ the
group of $\A$-linear automorphisms $f$ of $M$ satisfying
$\alpha \circ (f \times f) = \alpha$. As usually, if $M$ is the fibre
type of $\ell$, we say that {\em the structural
group of} $\ell$ {\em reduces to} $GL(M,\alpha)$, if $\ell$ has a
trivializing covering $\{ (U_i,\tau_i)\}_{i \in I}$, with
\begin{equation}
\tau_{jx} \circ \tau^{-1}_{ix} \in GL(M,\alpha),
\quad \forall i,j \in I, \; x \in U_i \cap U_j.
\end{equation}
If (2) holds, it is clear that the formula
\begin{equation}
s(x) := \alpha  \circ (\tau_{ix} \times \tau_{ix}).
\end{equation}
defines a hermitian form $(E_x,s(x))$ independent of the choice of
$i$, and that the induced map $s : E \oplus E \ra \A$ locally
coincides with $\alpha \circ (\tau_i \times \tau_i)$, hence it is an
$\A$-map. As a result, we have

\begin{The}
Let $\ell$ be an $\A$-bundle of fibre type $M$ and $(M,\alpha)$ a
hermitian form. If the structural group of $\ell$ reduces to
$GL(M,\alpha)$, then $\ell$ has an $\A$-hermitian structure.
\hfill $\Box$
\end{The}

If $\A$ is a complete C$^*$-algebra, the converse of Theorem 4.7 is
true in the topological case (\cite{MP;Red}). Besides, the reduction
of the structural group of a bundle does not depend on the choice of
$\alpha$ (ibid.).

\section{{\BBb A}-connections}

In this section, we define $\A$-connections as operators between the
sections of certain $\A$-bundles and we construct such a connection on 
a bundle having a trivializing covering with a subordinate
$\A$-partition of unity (Theorem 5.3). This connection is compatible
with the hermitian structure obtained in Theorem 4.7 (Theorem 5.4). 

\begin{Pro}
Let $(X,{\cal A}) \in Man_h(\A)$ modelled on $M$ and let $\ell =
(E,\pi,X)$ be an $\A$-bundle of fibre type $N$. If $\tilde{\ell} :=
(L(TX,E),\tilde{\pi},X)$, where
\[
L(TX,E) := \bigcup_{x \in X} L_{\A}(T(X,x),E_x)
\]
and $\tilde{\pi} : L(TX,E) \to X$ is the natural projection, then
$\tilde{\ell}$ admits the structure of an $\A$-bundle of fibre type
$P := L_{\A}(M,N)$.
\end{Pro}

\noindent {\em Proof.}
If $\{(U_i,\phi_i)\}_{i \in I}$ is an atlas of $X$ and
$\{(U_i,\tau_i)\}_{i \in I}$ a trivializing covering of $\ell$, let
\[
\tilde{\tau}_i : \tilde{\pi}^{-1}(U_i) \ra U_i \times P : f \mt
(\tilde{\pi}(f), \tau_{ix} \circ f \circ \bar{\phi}^{-1}_i),
\]
where $x = \tilde{\pi}(f)$ and $\bar{\phi}_i$ is given by (1). Then
$\{(U_i,\tilde{\tau}_i)\}_{i \in I}$ is a trivializing covering of
$\tilde{\ell}$. \hfill $\Box$

\begin{Def}
{\em Let $\ell = (E,\pi,X)$ be an $\A$-bundle. We denote by
$\Gamma(X,E)$ and $\Gamma(X,L(TX,E))$ the $\A$-modules of the
$\A$-differentiable sections of $\ell$ and of $\tilde{\ell}$,
respectively. We say that a mapping
\[
D : \Gamma(X,E) \ra \Gamma(X,L(TX,E)),
\]
is an} $\A$-connection on $\ell$, {\em if it is $\A$-linear and it
satisfies the} Leibniz condition:
\[
D{f\xi} = Tf \cdot \xi + f \cdot D\xi,
\]
{\em for any $\xi \in \Gamma(X,E)$ and any $\A$-map $f : X \ra \A$.

Besides, we say that $D$ is}  compatible {\em with a hermitian
structure $g$ of $\ell$, if 
\[
g_x(D\xi(x)(v),\eta(x)) + g_x(\xi(x),D_{\eta}(x)(v)) =
T_x(g \circ (\xi,\eta))(v),
\]
for every $x \in X$, $\xi, \eta \in \Gamma(X,E)$ and} $v \in T(X,x)$.
\end{Def}

\begin{The}
Let $\ell$ be an $\A$-bundle having a locally finite
trivializing covering with a subordinate $\A$-partition of unity.
Then $\ell$ has an $\A$-connection.
\end{The}

\noindent {\em Proof.}
(i) Assume first that the fibre type of $\ell$ is a free finitely
generated $\A$-module $\A^m$. Let $\{(U_i,\tau_i)\}_{i \in I}$ be the
locally finite trivializing covering of $\ell$ and
$\{(U_i,\psi_i\}_{i \in I}$ the subordinate $\A$-partition of unity.
For every $i \in I$, we set
\begin{equation}  \label{frame}
\varepsilon_{ij} : U_i \ra E_i := \pi^{-1}(U_i) : x \mapsto
\tau^{-1}_i(x,b_j)
\end{equation}
where $\{b_j\}_{1\leq j \leq m}$ is an arbitrary basis of $\A^m$. Then
$\{\varepsilon_{ij}\}_{1 \leq j \leq m}$  is an
$\A$-dif\-fer\-ent\-iable local frame of $\ell$ and every
$\xi \in \Gamma(U_i,E_i)$  is written as $\xi = \sum_j \xi_{ij}
\cdot \varepsilon_{ij}$, where $\xi_{ij} : U_i \ra \A$
is an $\A$-map, for every $j = 1,..., m$. We set
\[
D_i : \Gamma(U_i,E_i) \ra \Gamma(U_i,L(TU_i,E_i)) :
\xi \mapsto D_i(\xi) := \sum_j (T\xi_{ij}) \cdot \varepsilon_{ij}.
\]
It is straightforward that $D_i$ is a local $\A$-connection and that
\[
D : \Gamma(X,E) \ra \Gamma(X, L(TX,E)) : \xi \mapsto D\xi :=
\sum_i \psi_i \cdot D_i(\xi)
\]
is an $\A$-connection on $\ell$.

(ii) Suppose now that the fibre type of $\ell$ is $M \in \P$ and let
$N \in \P$ and $m \in \N$ with $M \oplus N = \A^m$. We consider the
trivial  $\A$-bundle $\ell_o := (X \times N, pr_1, X)$  and the
Whitney sum $\ell \oplus \ell_0 = (F,\bar{\pi},X)$. Let
$\ell \oplus \ell_o$ be endowed with the $\A$-connection
$\tilde{D} : \Gamma(X,F) \ra \Gamma(X,L(TX,F))$, obtained in (i).
If $I : E \ra F$ denotes the canonical injection and
$Pr : F \ra E$ the canonical projection, we obtain an $\A$-connection
$D$ on $\ell$, setting

\medskip
\hspace{2cm} $D\xi : X \ra L(TX,E) : x \mapsto D\xi(x) := Pr \circ
\tilde{D}(I \circ \xi)(x). \hfill \Box$

\medskip
If $\A$ is a complete $C^*$-algebra, the $\A$-hermitian structures
obtained in Theorem 4.7 and the $\A$-connections of the above
Theorem 5.3 are compatible. In fact, we have

\begin{The}
Let $\A$ be a commutative complete lmc C*-algebra with unit and $\ell$
an $\A$-bundle of fibre type $M$. Suppose that $\ell$ has a locally
finite trivializing covering $\{(U_i,\tau_i)\}_{i \in I}$
admitting a subordinate $\A$-partition of unity and satisfying (2),
for a hermitian form $(M,\alpha)$. Then $\ell$ has an $\A$-hermitian
structure and a compatible $\A$-connection.
\end{The}

\noindent {\em Proof.} 
(i) Assume first that the fibre type of $\ell$ is $\A^m$. Let
$\{x_j\}_{1 \leq j \leq m}$ be an orthonormal basis of $\A$ with
respect to $\alpha$ (cf. Cor. 4.2) and, for every $i \in I$, let
$\{\varepsilon_{ij}\}_{1 \leq j \leq m}$ be the
induced local frame of $\ell$ (see (4)). Let now $g$ be the 
$\A$-hermitian structure of $\ell$, constructed in Theorem 4.7, and
$D$ the $\A$-connection obtained in Theorem 5.3. Then, for every
$i \in I$,
\[
g_x(D_i\xi(x)(v),\eta(x)) + g_x(\xi(x),D_i\eta(x)(v)) =
\]
\[
= \sum_j (T_x \xi_{ij}(v) \cdot (\eta(x))^* +
       \xi_{ij}(x)(T_x\eta_{ij}(v)^*),
\]
where  $\xi|_{U_i} = \sum_j  \xi_{ij}e_{ij}$,  $\eta|_{U_i} =
\sum_j \eta_{ij}e_{ij}$, $x \in U_i$ and $v \in T(X,x)$. Thus,

\noindent
(5) \hfill \begin{minipage}{11cm}
\[
g_x(D\xi(x)(v),\eta(x)) + g_x(\xi(x),D\eta(x)(v)) =
\]
\[
= \sum_i \psi_i(x) \cdot \sum_j T_x(\xi_{ij} \cdot (\eta_{ij})^*)(v).
\]
\end{minipage}

\noindent
On the other hand,
\[
g \circ (\xi,\eta)(x) = g( \sum_j \xi_{ij}(x)\varepsilon_{ij}(x),
\sum_j \eta_{ij}(x)\varepsilon_{ij}(x))
= \sum_j \xi_{ij}(x) \cdot \eta_{ij}(x)^*,
\]
for every $(U_i,\tau_i)$ containing $x$, consequently,
\setcounter{equation}{5}
\begin{equation}
T_x(g \circ (\xi,\eta))(v) = \sum_j T_x(\xi_{ij} \cdot
         (\eta_{ij})^*)(v),
\end{equation}
implying the required equality.

(ii) Suppose now that the fibre type of $\ell$ is $M \in \P$. Let
$N \in \P$ and $m \in \N$ with $M \oplus N \equiv \A^m$. We consider
the trivial $\A$-bundle $\ell_o =
(X \times N,pr_1,X)$ and the sum  $\ell \oplus \ell_o$. The latter
is an $\A$-bundle of fibre type $\A^m$, the structural group of which
reduces to $GL(M \times N,\alpha\oplus\alpha_o)$, for any positive
definite $\A$-hermitian inner product $\alpha_o$ on $N$. Consider
$\ell$ and $\ell\oplus\ell_o$ provided with the $\A$-hermitian
structures $g$ and $\widetilde{g}$, where
\[
g(x) := \alpha  \circ  (\tau_{ix}  \times  \tau_{ix}) \ ; \;  i \in I,
\; x \in U_i,
\]
\[
\widetilde{g}(x) := g(x) \oplus \alpha_o \ ,  \; x \in X.
\]
Then, if $\{x_j\}_{1 \leq j \leq n}$ is again an orthonormal basis of
$\A$ with respect to $\alpha\oplus\alpha_o$, and, for every $i \in I$,
$(U_i, \Phi_i := \tau_i \times  id_N)$ is the trivializing
pair of $\ell \oplus \ell_o$ induced by $(U_i,\tau_i)$ of $\ell$, the
sections  $\varepsilon_{ij}(x) := \Phi^{-1}_i(x,x_j)$,
$(x \in U_i, j = 1,...,n)$ form an orthonormal $\A$-differentiable
frame on $U_i$. As in (i), we construct an $\A$-connection
$\widetilde{D}$ on $\ell \oplus \ell_o$ which is compatible with
$\widetilde{g}$. Now, setting
\[
D : \Gamma(X,E) \ra \Gamma(X,L(TX,E)) : \xi \mapsto
Pr(\widetilde{D}(I \circ \xi))
\]
(see the proof of Theorem 3.3), we obtain
\begin{eqnarray*}
g_x(D\xi(x)(v),\eta(x)) & + & g_x(\xi(x),D\eta(x)(v)) = \\
         & = & g_x(Pr \circ \widetilde{D}(I \circ \xi)(x)(v),
                   I \circ \eta(x)) + \\
         &   & \quad + g_x(I \circ \xi(x),
                   Pr \circ \widetilde{D}(I \circ \eta)(x)(v)) = \\
         & = & g_x(Pr \circ \widetilde{D}(I \circ \xi)(x)(v),\eta(x))
                + \\
         &   & \quad + \alpha_o((1-Pr) \circ \widetilde{D}(I
                     \circ \xi)(x)(v),0) + \\
         &   & \quad + g_x(\xi(x), Pr \circ \widetilde{D}(I
                    \circ \eta)(x)(v)) + \\
         &   & \quad + \alpha_o(0, (1-Pr) \circ (\widetilde{D}(I
                     \circ \eta)(x)(v)) = \\
         & = & \widetilde{g}_x(\widetilde{D}(I \circ \xi)(x)(v),I
                      \circ \eta(x)) + \\
         &   & \quad + \widetilde{g}_x(I \circ \xi(x),
                     \widetilde{D}(I \circ \eta)(x)(v)) = \\
         & = & T_x(\widetilde{g} \circ (I \circ \xi,
                     I \circ \eta))(v) = \\
         & = & T_x(g \circ (\xi, \eta))(v),
\end{eqnarray*}
which completes the proof. \hfill $\Box$

\noindent
Department of Mathematics \\
University of Athens \\
Panepistimiopolis \\
Athens 157 84, Greece \\
e-mail:mpapatr@atlas.uoa.gr

\end{document}